\documentclass[conference]{IEEEtran}
\usepackage[]{times}
\usepackage{epsfig}
\usepackage{subfigure}
\usepackage{amsmath}
\usepackage{bm,amssymb}
\usepackage{graphicx}
\usepackage[top=0.5in, bottom=0.75in, left=0.75in, right=0.75in]{geometry}
\usepackage{multirow}
\usepackage{flushend}
\usepackage{url}
\usepackage{color}

\title{Robust Broadcast-Communication Control of Electric Vehicle Charging}
\author{
\authorblockN{Konstantin Turitsyn}
\authorblockA{CNLS \& Theoretical Divison\\
Los Alamos National Lab\\
Los Alamos,\\ NM 87545, USA\\
Email: turitsyn@lanl.gov}
\and
\authorblockN{Nikolai Sinitsyn}
\authorblockA{Theoretical Divison\\
Los Alamos National Lab\\
Los Alamos,\\ NM 87545, USA\\
Also with NMC\\
Email: nsinitsyn@lanl.gov}
\and
\authorblockN{Scott Backhaus}
\authorblockA{MPA Divison\\
Los Alamos National Lab\\
Los Alamos,\\ NM 87545, USA\\
Email: backhaus@lanl.gov}
\and
\authorblockN{Michael Chertkov}
\authorblockA{CNLS \& Theoretical Divison\\
Los Alamos National Lab\\
Los Alamos,\\ NM 87545, USA\\
Also with NMC\\
Email: chertkov@lanl.gov}}

\begin{document}
\maketitle
\begin{abstract}
The anticipated increase in the number of plug-in electric vehicles (EV) will put
additional strain on electrical distribution circuits.  Many control schemes have been
proposed to control  EV charging.  Here, we develop control algorithms based on
randomized EV charging start times and simple one-way broadcast communication allowing for a time delay between communication events.  Using arguments from queuing theory and
statistical analysis, we seek to maximize the utilization of excess distribution circuit capacity while keeping the probability of a circuit overload negligible.  


\end{abstract}

%\IEEEpeerreviewmaketitle

\section{Introduction}
\label{sec:intro}

If electric vehicles (EV) make up a significant fraction of the future vehicle fleet,
there will be significant impacts on electrical generation and transmission due to the
increase in electricity consumption.  Some researchers have viewed the potential
increase in EVs as an opportunity to utilize the on-board battery storage in an
interactive way to provide two-way energy flows to buffer time-variable renewables\cite{09CF,97KL,05KT} or to provide ancillary services such as (short-term) frequency regulation\cite{02V2G}.  We also view
these potential applications as mutually beneficial to consumers, utilities, and
generators, however, the first potential impact of EVs,
independent of mass deployment of renewable generation, is filling in the nighttime
load-curve minimum typical of most regions.\cite{07PNNL} The electrical system as it is currently configured would benefit in several ways.  The primary benefit is a higher
utilization of currently deployed assets resulting in faster paybacks on generation,
transmission, and distribution investments. Second, in some regions, the locational
marginal prices for wholesale energy become negative during the nighttime hours.  If
properly controlled, EVs will provide a reliable increased nighttime load making the
operation of baseload generation more profitable.  Modifying the load curve in this way
will also reduce the relative impact of demand charges relative to energy charges for
many distribution utilities.

To serve this new load, distribution utilities must determine how many EVs can be
reliably served on a given distribution circuit without substantially increasing the
probability of an overload or other system upset. A likely overload scenario is the
synchronization of EV charging start times that would occur as their owners begin EV
charging immediately after they arrive home from work during the weekday evening hours.  We seek to avoid this load synchronization by controlling EV charging start times.

Although sophisticated control schemes relying on high-speed two-way communication
between all EV loads and centralized controllers could allow distribution circuits to run at nearly full capacity, the reliability of electrical service would then
compromised due to real-time reliance upon a sophisticated communications network.  Instead, we seek a distributed control method that relies upon no or minimal communication, where the
decision of when to begin charging is made local to each EV charger.  We propose to
regulate EV charging start times by modulating a single, circuit-wide EV connection rate (or arrival rate) that determines, on average, how many EV chargers can commence charging per unit time. The circuit-wide rate may be sent to each charger via one-way, broadcast
communication and updated periodically at a frequency that suits a particular utility's or
circuit's needs.  As a limit of this scheme, we also consider a simpler version where a
fixed arrival rate is preprogrammed into the charger by the utility.

\begin{figure}
\centering \caption{a) Example schematic load curve showing a load valley lasting for time $\Delta T$ during the nighttime hours where excess capacity on the circuit can accommodate a maximum of $N$ EV chargers without overloading where each EV charger consumes $P_{con}$ real power. b) Schematic representation of control scheme.  The open circles represent the sparse measurements of $n_i$ at the beginning of a time interval $\tau$.  See main text for details.
} \includegraphics[width=0.55\textwidth]{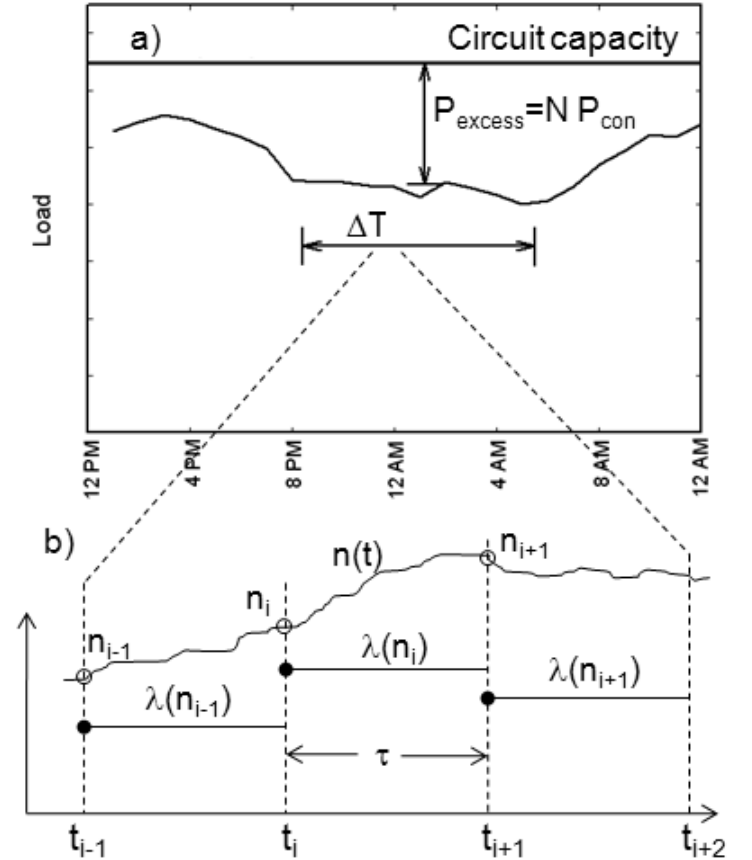}
\label{fig:schematic}
\end{figure}

For this initial work, we consider a few simplifying assumptions.  We envision that each utility or circuit will have a predetermined EV-charging `window' of length $\Delta T$
with start and stop times that are correlated with the minimum in the load curve for
that utility or circuit, schematically shown in Fig~\ref{fig:schematic}a.  We also
assume that, during $\Delta T$, the other loads are known well enough that an available
excess capacity $P_{excess}$ can be reliably determined and that all EV chargers consume the same amount of power $P_{con}$ so that the capacity for EV charging is expressed in terms of the maximum number of EVs charging without an overload, i.e. $N=P_{excess}/P_{con}$.  Each of these assumption can relaxed to allow for charging at
any point in the load curve,  fluctuations in uncontrolled, non-EV loads on the circuit, and multiple classes of EV charging powers.

Schematically shown in Fig~\ref{fig:schematic}b, the window $\Delta T$ is split into
many smaller time intervals of length $\tau$ with beginning and end times $t_i$.  At
each $t_i$, the number $n_i$ of EVs charging is measured{\footnote {To develop intuition and to keep the computations simple in this initial work, we consider the situation
where $n_i$ is available to the entity controlling the system.  Although this may
require two-way communication, this can easily be relaxed by instead measuring the real
or apparent power drawn by the circuit using a current transformer at the substation and controlling EV charging based the apparent power available below the circuit capacity.}} and $n_i$ is used to determine the connection (or arrival) rate of EV chargers
$\lambda(n_i)$ that will apply until the end of the internal at $t_{i+1}$, i.e. for a
time $\tau$.  In this context, $\tau$ is the time interval between utility broadcast
communication events.  During this interval, EVs are connecting at a rate $\lambda(n_i)$ but are also leaving at an average rate $\mu n(t)$ where $\mu=\Delta E/P_{con}$ is the
inverse of the average EV charging time and $\Delta E$ is the average amount of energy
required to reach full EV battery charge.  Although the individual decisions to start
and stop charging are made locally, control of the charging load is exercised through
broadcasting a single parameter $\lambda(n_i)$ which amounts to feedback control because it depends on the number of EV charging at time $t_i$.

The EV service concept proposed above is extremely robust.  The connection of EV
chargers is inherently a random process naturally avoiding any of the dangerous
synchronization that can occur in EV charging.  The reliance on the communication network is minimal since decisions to start or stop charging are distributed and based upon a single number $\lambda (n_i)$ received via rate-limited, one-way, broadcast communication. In addition, misoperation by one or a handful of EV chargers will not greatly affect the overall load.  However, the randomness of the EV connection process will cause the
aggregate EV load fluctuate (in addition to the fluctuations in the underlying load). Excluding the possibility of real-time intervention
during the time interval $\tau$, a utility must carefully choose $\lambda(n_i)$ to
maximize EV service and asset utilization without overloading or damaging the circuit.
The remainder of this manuscript analyzes this problem and seeks to make good choices for $\lambda(n_i)$.

\section{Calculations}
\label{sec:setup}
Within any period $\tau$, the EV charging process described above closely resembles the
birth-death queuing models typically used in telecommunications where both the arrival
and completion of calls is often assumed to be a random process with an average arrival rate $\lambda$ and departure (or end of service) rate $\mu$.  For our EV-charging
problem, one implication of adopting this model is that the time between the start of
EV-charing events is exponentially distributed.  The same holds for the time between
completion of EV-charging events.  In practice, the arrival and departure-time
distributions will certainly be cut off by $\Delta T$ and $\Delta E_0/P_{con}$,
respectively, where $\Delta E_0$ is the EV battery capacity.  However, for a large
number of EVs in the system, the effect of eliminating these tails should be minimal.

We begin with a simplified `constant-rate' version of the
problem letting $\tau \rightarrow \infty$ corresponding to the case where
$\lambda=\lambda_0$ is preprogrammed and the problem is essentially one of design and
not control.  In this limit, the problem reduces to an $M/M/\infty$ queuing
problem\cite{kleinrock} where the probability of finding $n$ EVs charging is given by
\begin{equation}\label{pk}
\pi_n=e^{-\bar{n}} \bar{n}^n/n!,
\end{equation}
and the probability of a circuit overload, i.e. $n\geq N$, is given by
\begin{equation}\label{overload1}
P_N=\sum_{m=N}^\infty \pi_m=1-\Gamma(N,\bar{n})/\Gamma(N).
\end{equation}
Here, $\bar{n}=\sum_n n \pi_n=\lambda_0/\mu$ is the average number of EVs charging and
$\Gamma (\cdot)$ and $\Gamma (\cdot,\cdot)$ are the Gamma and incomplete Gamma
functions.

Next, we consider a `variable-rate' version of the problem where we control EV charging
by allowing $\lambda$ to depend on $n_i$.  The goal of this control protocol is to
maximize the utilization of excess circuit capacity while keeping the probability of
overload $P_N$ acceptably low. For a given control protocol $\lambda(n)$, the probability distribution $\pi_n^i$ of the number $n_i=n$ of EVs charging at time $t_i$, obeys a simple Markov chain equation:
\begin{equation}
 \pi_n^i = \sum_m W_{nm}(\lambda(m)\tau) \pi_m^{i-1}
\end{equation}
where we have introduced $W_{mn}(\lambda\tau)$ - the transition kernel of the
$M/M/\infty$ process with the rates $\lambda$ and $\mu$.  For $n>m$, the transition kernel can be expressed as the following sum \cite{schoutens2000stochastic}:
\begin{equation}
   W_{nm}(\lambda\tau) = \exp(-q)\sum_{k=0}^n \frac{n! q^{n-k} \epsilon^{m-k}(1-\epsilon)^k}{k!(n-k)!(m-k)!},
\end{equation}
where $\epsilon = 1 - \exp(-\mu \tau)$ and $q = \lambda \epsilon/\mu$.
Assuming that the system was not overloaded at time $t_{i-1}$, the probability of an overload at time $t_i$ (i.e. $n_i \geq N$) is given by $P_N =
\sum_{m\geq N} \pi_m^i$ and can be written as
\begin{equation} \label{PNineq}
 P_N = \sum_{n<N,m\geq N} W_{mn}(\lambda(n)\tau) \pi_n^{i-1} \leq \max_{n<N} \sum_{m\geq
 N}W_{mn}(\lambda(n)\tau)
\end{equation}
The last inequality, which is based on the normalization of $\pi_n^{i-1}$ and that $\pi_n^{i-1}=0$ for $n \geq N$,  forms the foundation for constructing our variable-rate protocol. We define the control protocol $\lambda(n)$ as a solution of
\begin{equation}\label{lambdaeq}
\sum_{m\geq N} W_{mn}(\lambda(n)\tau)= P_N^* \ll 1, \quad n<N
\end{equation}
and set $\lambda(n)=0$ for $n\geq N$. For a given value of $P_N^*$ this protocol is
guaranteed to have an overload probability bounded from above: $P_N < P_N^*$. We note
that this control scheme might be suboptimal because we impose more strict constraints on the overload probability than required; instead of requiring the total overload probability to be less then $P_N^*$, we look for a solution where the probability of overload conditioned on the starting with $n_{i-1}$ EVs connected is bounded by $P_N^*$ for {\it every} $n_{i-1}$.
However, as we will show, the performance of this control scheme is almost perfect, so
further optimization can result only in diminishingly small improvements.  The result also indicates that one starting condition $n_{i-1}$ dominates the transition to the overloaded state.  

As a first measure of the control system's performance, we consider the time response of the control or equivalently the time to relax to a steady state.  In the constant-rate scheme, this time is estimated as the product of the average number of EVs charging in equilibrium and the typical inter-arrival time $\bar{n}/\lambda_0$, which according to the discussion below Eq.~(\ref{overload1}) is $\bar{n}/\lambda_0=1/\mu$.  Control in this case can be slow because the relaxation time is the typical time required to charge an EV, which can be several hours for present EV battery and charging technology\cite{07PNNL,10SAE}.  In the variable-rate case, the relaxation time can be estimated as the number of time intervals $\tau$ it takes, starting with $n_0=0$, to reach $n\sim N$. To leading order, the solution of (\ref{lambdaeq}) can be estimated as
$\lambda(n)\tau\approx N-n - \Delta_n$, where $\Delta_n \sim -\log P_N \sqrt{N-n}$ for
$N-n\gg 1$. Thus, in the limit $N\gg 1, \mu\tau \ll 1$ and $\mu N \tau \lesssim 1$,
there will be $n\approx N+\log P_N \sqrt{N}$ cars in the system after one interval $\tau$.  After
the second interval $\tau$, $n\approx N+ \log P_N \sqrt{-\log P_n \sqrt{P_N}}$. Thus, the total time
required to reach stationary state can be estimated as $\tau \log (-N \log P_N)$.
For realistic values of $N$ and $P_N$, the response time of the variable-rate protocol is now under the control of the utility which can set $\tau$ based on technological or economic considerations.

A second measure of system's performance, which we address in the remainder of this
manuscript, is the average capacity utilization defined as $\bar{n}/N$ ranging from $0$
with no EVs charging to $1$ when the excess circuit capacity is fully utilized.  Our
control schemes seek to maximize $\bar{n}/N$ while keeping $P_N$ acceptably small.  The
proximity of $\bar{n}/N$ to one for a given $P_N$ is a measure of the quality of the
control scheme.  Exact calculation of this quantity is a computationally challenging
problem that is beyond the scope of this work. Instead, here we obtain the lower bound
for the capacity by exploiting the convexity of the function $\lambda(n)$ derived from Eq.~(\ref{lambdaeq}). We first note, that the stability of the queue (i.e.
existence of the stationary distribution) implies that the following relation holds:
\begin{equation}\label{conslaw}
 \overline{\lambda(n)\tau} = \mu\tau \bar{n}
\end{equation}
It is not possible to extract the value of $\bar{n}$ from Eq.~(\ref{conslaw}) directly because $\lambda(n)$ is a nonlinear function of $n$. However, for convex and monotonically
decreasing functions $\lambda(n)$ that solve the equation (\ref{lambdaeq}) one can use
the relation (\ref{conslaw}) to find the lower bound for $\bar{n}$. It follows from
Jensen's inequality that $\mu \bar{n}-\overline{\lambda(n)} \geq \mu
\bar{n}-\lambda(\bar{n})$ and from motonicity of $\mu n - \lambda(n)$ we conclude that
$\bar{n}\geq n^*$ with $n^*$ being the solution of the equation
\begin{equation}\label{nbareq}
 \lambda(n^*) =\mu n^*,
\end{equation}
where we have assumed some convex and monotonic continuation of the discrete function
$\lambda(n)$.

Equations~(\ref{lambdaeq}) and (\ref{nbareq}) have been solved numerically with
Newton's method with the help of the approximate expressions for the transition kernel
$\sum_{m\geq N}W_{mn}(\lambda \tau)$. For a given value of $P_N^*$, a table of
$\lambda(n)$ was found for all $n=0\dots N-1$ which was then continued to
real values of $n$ with second-order spline interpolation. Equation~(\ref{nbareq})
was then solved to find the lower bound on $\bar{n}$, i.e. $n^*$, as a function of the
overload probability bound $P_N^*$.

\section{Results and Discussion}
\label{sec:results}
\begin{figure}
\centering \caption{EV connection rate $\lambda(n)$ versus $n$ for the constant-rate (red, dashed line) and variable-rate cases (blue, solid curve) for $N=100$ and $P_N = 10^{-10}$.  For the variable-rate case, $\mu \tau=0.001$.
} \includegraphics[width=0.49\textwidth]{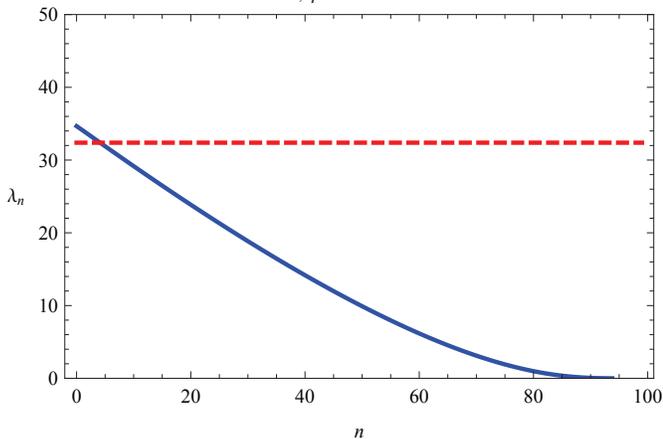}
\label{fig:protocols}
\end{figure}

To understand the differences between the constant-rate and variable-rate protocols, Fig.~\ref{fig:protocols} displays $\lambda(n)$ versus $n$ for $N=100$ and $P_N=10^{-10}$.  For the variable-rate case, the condition on the overload probability in Eq.~(\ref{lambdaeq}) suppresses $\lambda(n)$ as $n$ approaches $N$ while no such condition exists in the constant-rate case.  The resetting of $\lambda(n)$ at the beginning of each time interval $\tau$ provides a mechanism for controlling the fluctuations in $n(t)$ due to the random EV connection process.  For example, if fluctuations in the rate of EVs connections between $t_{i-1}$ and $t_i$ generated a few more EV connections that was expected on average, i.e. more than $\tau \lambda(n_{i-1})$, then average connection rate in the subsequent period is decreased.  The opposite control response occurs if the average number of connections falls short of the expected average.

We find that the particular shape of $\lambda(n)$ as $n$ approaches $N$ is important.  In Fig.~\ref{fig:protocols}, $N=100$ and $\lambda(n)$ could be closely approximated by a linear function up to about $n=60$.  Extending that line to $\lambda=0$ near $n=70$ would provide a very familiar type of control, i.e. proportional feedback control with a deadband between $n=70$ and $n=100$.  However, our analysis of this modified control protocol finds this linearized $\lambda(n)$ shows little improvement in average capacity utilization over the constant-rate case demonstrating that the details of $\lambda(n)$ for $n\sim N$ are crucial to optimal control performance.

\begin{figure}
\centering \caption{Average capacity utilization for the constant-rate case (red, dash-dot curve) and the variable-rate case with $N=100$ and $\mu \tau = 0.01$ (blue, dashed curve) and $\mu \tau = 0.001$ (blue, solid curve).  The horizontal axis is the common log of the inverse of the probability of overload. Intuitively for all cases, the average capacity utilization falls as the probability of overload become smaller.  In addition, the utilization also falls as $\tau$ increases in the variable-rate case.
} \includegraphics[width=0.49\textwidth]{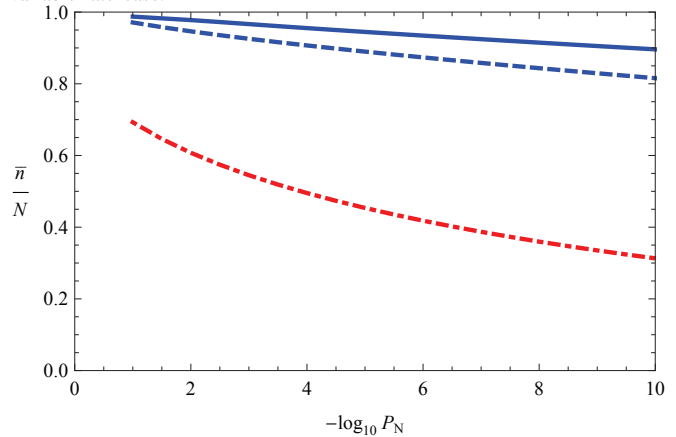}
\label{fig:plt100}
\end{figure}

Figure~\ref{fig:plt100} shows the average capacity utilization for both the constant and variable-rate cases.  The ability of the variable-rate control to essentially clip fluctuations in $n$ by modifying $\lambda$ is clearly evident in the the higher values of $\bar{n}/N$ achieved.  Figure~\ref{fig:plt100} also demonstrates the incremental value of broadcast-communication control over the case of no communication.  Considering low overload probabilities, e.g. $P_N=10^{-10}$, Fig.~\ref{fig:plt100} shows even when the optimal value of $\lambda_0$ is used, the best excess capacity utilization in the constant-rate scheme is only $\approx 30\%$.  Simply adding one-way broadcast communication to allow the utility to update the connection rate every $\tau=0.01/\mu$, the utilization jumps to $\approx 80\%$.  To provide a sense of the communication time scale, a typical charging time for an EV using today's technology\cite{10SAE} is approximately 5 hours.  Therefore, to achieve $80\%$ utilization at $P_N=10^{-10}$ would require a communication rate of once every 3 minutes, which is achievable even with slow communication methods such as power line carrier (PLC).  Figure~\ref{fig:plt100} shows that increasing the communication rate by a factor of ten (about once every 20 seconds) only boosts the average capacity utilization from $80\%$ to $90\%$.   From this analysis, we conclude that the presence of relatively slow, one-way, broadcast communication is sufficient to allow utilities to effectively manage EV charging.

\begin{figure}
\centering \caption{Same as Fig.~\ref{fig:plt100} except here $N=1000$.  For larger $N$ using variable-rate control, the capacity utilization approaches 1 and becomes less sensitive to $\tau$.
} \includegraphics[width=0.49\textwidth]{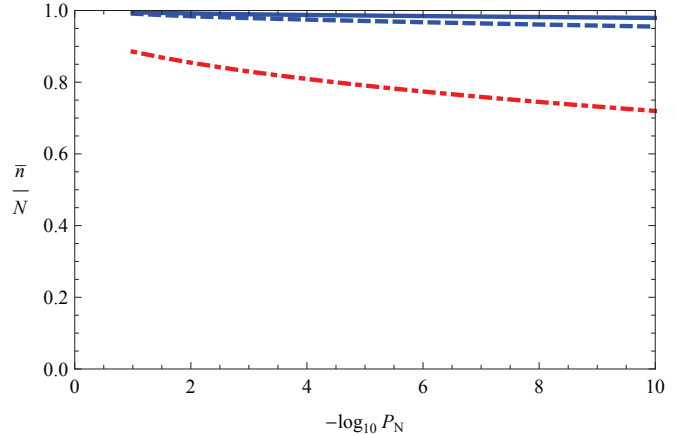}
\label{fig:plt1000}
\end{figure}

Figure~\ref{fig:plt1000} provides the same analysis as Fig.~\ref{fig:plt100} except the available excess capacity $N$ has been increased from 100 to 1000 EVs.  The performance of both control methods improves as $N$ increases.  In each case, the better performance is due to fractionally smaller variance in the equilibrium distribution $\pi_n$ for the constant-rate case and in the number of EV connections during $\tau$ for the variable-rate.  As $N$ is increased, the value of even slow communication is still apparent due to the increase in excess capacity utilization from $70\%$ to about $95\%$. However, increasing the speed of communication becomes less and less important.

\section{Conclusion and Future Directions}
\label{sec:conclusions}
Avoiding synchronization of EV-charging start times is crucial to avoiding circuit overloads as the penetration of EVs on a distribution circuit increases.  Randomization of start times is one method that avoids this synchronization and also allows the properties of the EV charging load to be analyzed from the standpoint of queuing  theory.  Our initial analysis shows that even relatively slow, one-way, broadcast communication such as PLC can be used to modulate the random EV connection process and effectively control the EV-charging load.  When there is a large amount of excess circuit capacity for EV charging (e.g. 1000 EV connected at a time), our control algorithm easily achieves $95\%$ utilization of the capacity with a minimal probability of overloading the circuit.  Our analysis also shows that compared to no communication, even relatively slow, one-way, broadcast communication boosts the utilization of circuit capacity dramatically.  However, as the speed of that communication increases, there are diminishing returns for EV charging.

Our initial work lays the foundation for further study into this EV-charging control method that could also be applied to other forms of demand response if the consumer finds randomization of load start times acceptable.  Our work can be expanded on in several ways:
\begin{itemize}
\item{To avoid having to define an accurate excess circuit capacity for every distribution circuit, a utility may choose to use measurements of the current injected into the circuit at the substation to determine excess capacity in real time.  These measurements would also capture the statistical fluctuations in the non-EV-charging loads.  Our analysis should be performed in this setting to account for time-dependent average excess capacity and its fluctuations.}

\item{Our analysis assumed that the population of EVs connected to their chargers but not presently charging was relatively large so that the number of unconnected EVs had little dependence on the number of connected EVs.  During daylight hours when EVs are driven to the workplace, the number of EVs on a residential distribution circuit may be greatly reduced and also fluctuate greatly as they are used for day trips.  Our analysis should be extended to account for these effects.}
\item{Our analysis assumes that once an EV begins charging, it stays connected until its battery is fully charged.  A more realistic setting would include the possibility of EVs disconnecting before reaching full charge either due to the requirements of the owner or perhaps due to a broadcast control signal from the utility.}
\item{The constraint in our model is based only on circuit capacity and does not include any effects of power flow.  Including the effects of reactive power would make for a much richer problem because EVs at different locations on the circuit would contribute different amounts of current draw at the substation.}
\item{Our model was confined to a single distribution circuit.  Natural extensions would consider both capacity and power flow effects at substation or even local transmission scales.}
\item{Finally, a more general model would also incorporate the transportation grid and may account for delays in transportation flows,  EVs changing stations (e.g. from home to work and back) and possibly even cross-correlations between irregularities in traffic flows and EV charging.}
\end{itemize}

\section*{Acknowledgment}

We are thankful to all the participants of the ``Optimization and Control for Smart
Grids" LDRD DR project at Los Alamos
and Smart Grid Seminar Series at CNLS/LANL for multiple fruitful discussions.  Research
at LANL was carried out under the auspices of the National Nuclear Security
Administration of the U.S. Department of Energy at Los Alamos National Laboratory under
Contract No. DE C52-06NA25396.

\bibliographystyle{IEEEtran}
\bibliography{SmartGrid}

\end{document}